\documentclass[reqno]{amsart}
\usepackage{amssymb}
\usepackage{epsf}

\newtheorem{theorem}{Theorem}

\theoremstyle{definition}

\theoremstyle{remark}
\newtheorem{remark}{Remark}

\begin{document}

\title[Unit circle elliptic beta integrals]
{Unit circle elliptic beta integrals}

\author{J.F. van Diejen}
 \address{Instituto de Matem\'atica y F\'{\i}sica,
 Universidad de Talca, Casilla 747, Talca, Chile}

\author{V.P. Spiridonov}
 \address{Bogoliubov Laboratory of Theoretical Physics, Joint
Institute for Nuclear Research, Dubna, Moscow Region 141980, Russia}

\thanks{{\em Date:} August 2003; {\em Ramanujan J.}, to appear}

\begin{abstract}
We present some elliptic beta integrals with a base parameter on
the unit circle, together with their basic degenerations.
\end{abstract}

\maketitle

\section{Introduction}
The theory of generalized gamma functions has been set up by
Barnes \cite{bar:multiple}. A slightly different approach was
advocated by Jackson, who considered the basic gamma function
depending on one base parameter $q$ and the elliptic gamma
function depending (symmetrically) on two bases $p$ and $q$
\cite{jac:basic}. For a long time, only the first of these
generalized gamma functions was appreciated in the literature
\cite{aar:special}. Recently, however, the elliptic gamma
function also got appropriate attention after the work of
Ruijsenaars \cite{rui:first}, who introduced it in the context
of integrable systems and investigated some of its properties. A
further study of the function in question was conducted by Felder
and Varchenko \cite{fel-var:elliptic}.  A modified elliptic gamma
function, which admits analytic continuation in one of the base
parameters, e.g. $q$, onto the unit circle $|q|=1$, has been
introduced recently by one of us in \cite{spi:integrals}.

In this paper we study beta type integrals on the unit circle
$|q|=1$ built of modified elliptic gamma functions, as well as
their basic degenerations. The first exact beta
type integration formula involving the conventional elliptic gamma
function was discovered in \cite{spi:elliptic}. Various
multidimensional generalizations of this elliptic beta integral
associated with the $C_N$ and $A_N$ root systems have been
investigated in \cite{die-spi:elliptic,die-spi:selberg,
rai:privite,spi:integrals}. A general theory of theta hypergeometric
integrals on tori and the Jacobi theta function generalizations of
the Meijer function was developed in \cite{spi:integrals}. The
beta integrals considered below should be thought of as the
$|q|=1$ counterparts of the integrals in \cite{spi:elliptic} and
\cite{die-spi:elliptic}. Recently, the Bailey's technique of deriving
identities for series of hypergeometric type \cite{aar:special}
has been generalized to the level of integrals \cite{spi:bailey}.
It can also be extended to the $|q|=1$ integrals under discussion.

\section{Preliminaries:\\ the Jacobi theta function and the elliptic
gamma function}

The main underlying structural object of this
paper is a Jacobi type theta function defined as
\begin{equation}\label{theta}
\theta(z;p)=(z;p)_\infty (pz^{-1};p)_\infty, \quad
(a;p)_\infty=\prod_{n=0}^\infty(1-ap^n) ,
\end{equation}
with $z,p\in\mathbb{C},|p|<1.$ It satisfies the transformation
properties
\begin{equation}\label{fun-rel}
\theta(pz;p)=\theta(z^{-1};p)=-z^{-1}\theta(z;p).
\end{equation}
Evidently, $\theta(z;p)=0$ for $z=p^{m},\, m\in\mathbb{Z}$, and
$\theta(z;0)=1-z$.  If we denote $p=e^{2\pi i\tau}$,
$\text{Im}(\tau)>0$, then the standard Jacobi $\theta_1$-function
\cite{whi-wat:course} is related to $\theta(z;p)$ as
\begin{eqnarray} \nonumber
\theta_1(u|\tau) &=&
-i\sum_{n=-\infty}^\infty (-1)^np^{(2n+1)^2/8}e^{\pi i(2n+1)u}
\\
&=& p^{1/8} ie^{-\pi iu}\: (p;p)_\infty\: \theta(e^{2\pi iu};p),
 \quad  u\in\mathbb{C}.
\label{theta1}\end{eqnarray} The complex function
$\theta_1(u|\tau)$ is entire, odd, and doubly quasiperiodic in $u$
\begin{eqnarray}\nonumber
&& \theta_1(u+1|\tau)     = -\theta_1(u|\tau),   \\
&& \theta_1(u+\tau|\tau) = -e^{-\pi i\tau-2\pi iu}
\theta_1(u|\tau). \label{quasi}\end{eqnarray} The modular
$PSL(2,\mathbb{Z})$-group,
\begin{equation}\label{sl2}
\tau\to \frac{a\tau+b}{c\tau+d},\qquad u\to \frac{u}{c\tau+d},
\end{equation}
with $a,b,c,d \in\mathbb{Z}$ and $ad-bc=1$, is generated by the
two transformations $\tau\to\tau+1$, $u\to u$ and $\tau\to
-\tau^{-1}$, $u\to u\tau^{-1}$. Its action on the Jacobi theta
function is determined by
\begin{subequations}
\begin{eqnarray}\label{modular-tr1}
&& \theta_1(u|\tau+1) = e^{\pi i/4}\theta_1(u|\tau),\qquad
\\  \label{modular-tr2}
&& \theta_1\big(\frac{u}{\tau}|-\frac{1}{\tau}\big)
=-i(-i\tau)^{1/2}e^{\pi iu^2/\tau} \theta_1(u|\tau) .
\end{eqnarray}
\end{subequations}
(Throughout this paper
the sign of the square root is fixed in accordance with the
principal branch with the cut chosen on the negative real axis.)
{}From the second of these relations, combined with the modular
transformation law for the Dedekind $\eta$-function
\begin{equation}
e^{-\frac{\pi i}{12\tau}} \left(e^{-2\pi i/\tau};e^{-2\pi
i/\tau}\right)_\infty =(-i\tau)^{1/2}e^{\frac{\pi i\tau}{12}}
\left(e^{2\pi i\tau};e^{2\pi i\tau}\right)_\infty ,
\label{ded}\end{equation} one readily deduces a corresponding
modular transformation formula for the $\theta(z;p)$ function
\begin{equation}
\theta(e^{2\pi i\frac{u}{\tau}};e^{-2\pi i\frac{1}{\tau}})
=-ie^{\pi i(\frac{u^2}{\tau}+\frac{\tau}{6}+\frac{1}{6\tau}
+\frac{u}{\tau}-u)} \theta(e^{2\pi iu};e^{2\pi i\tau}).
\label{mod}\end{equation}

The elliptic gamma function $\Gamma(z;q,p)$, $|q|, |p|<1$, is
related to the above theta function through the difference
equations
\begin{equation}
\Gamma(qz;q,p)=\theta(z;p)\Gamma(z;q,p),\qquad
\Gamma(pz;q,p)=\theta(z;q)\Gamma(z;q,p).
\end{equation}
It is given by the explicit product representation \cite{rui:first}
\begin{equation}
\Gamma(z;q,p) = \prod_{j,k=0}^\infty\frac{1-z^{-1}q^{j+1}p^{k+1}}
{1-zq^jp^k}. \label{ell-gamma}\end{equation} From this
representation the following reflection relation is immediate
\begin{equation}\label{refl}
\Gamma(z;q,p)\Gamma(z^{-1};q,p)=\frac{1}{\theta(z;p)\theta(z^{-1};q)}.
\end{equation}

\section{The modified elliptic gamma function}

The modified elliptic gamma function introduced in
\cite{spi:integrals} is constructed as a product of two elliptic
gamma functions of the form in (\ref{ell-gamma}), corresponding to
two different choices of bases. It is convenient to pass to an
additive formulation by introducing three pairwise incommensurate
quasiperiods $\omega_1$, $\omega_2$, $\omega_3$ and write
\begin{eqnarray}\nonumber
&& q=e^{2\pi i\frac{\omega_1}{\omega_2}}, \qquad p=e^{2\pi
i\frac{\omega_3}{\omega_2}},\qquad r=e^{2\pi
i\frac{\omega_3}{\omega_1}},
\\
&& \tilde q =e^{-2\pi i\frac{\omega_2}{\omega_1}}, \qquad \tilde
p=e^{-2\pi i\frac{\omega_2}{\omega_3}},\qquad \tilde r=e^{-2\pi
i\frac{\omega_1}{\omega_3}}, \label{ell-bases}\end{eqnarray} (i.e.,
$\tau=\omega_3/\omega_2$). The tilded bases $\tilde q, \tilde p,$
and $\tilde r$ are the respective modular transformations of
$q,p,$ and $r$. For
$\text{Im}(\omega_1/\omega_2),\text{Im}(\omega_3/\omega_1),\text{Im}(\omega_3/\omega_2)>0$
(so $|q|, |p|, |r|<1$), the modified elliptic gamma function is
now defined as \cite{spi:integrals}
\begin{equation}
G(u;\boldsymbol{\omega})= \prod_{j,k=0}^\infty \frac{(1-e^{-2\pi
i\frac{u}{\omega_2}}q^{j+1}p^{k+1}) (1-e^{2\pi i
\frac{u}{\omega_1}}{\tilde q}^{j+1}{r }^k)} {(1-e^{2\pi i
\frac{u}{\omega_2}}q^jp^k) (1-e^{-2\pi i
\frac{u}{\omega_1}}{\tilde q}^j{r}^{k+1})}.
\label{ell-d}\end{equation} It satisfies three difference
equations
\begin{subequations}
\begin{eqnarray}
&& G(u+\omega_1;\boldsymbol{\omega})=\theta(e^{2\pi
i\frac{u}{\omega_2}};p) G(u;\boldsymbol{\omega}),
\label{ell-1eq} \\
&& G(u+\omega_2;\boldsymbol{\omega})=\theta(e^{2\pi
i\frac{u}{\omega_1}};r) G(u;\boldsymbol{\omega}),
\label{ell-2eq} \\
&& G(u+\omega_3;\boldsymbol{\omega})= \frac{\theta(e^{2\pi
i\frac{u}{\omega_2}};q)} {\theta(e^{2\pi
i\frac{u}{\omega_1}}\tilde q;\tilde q)} G(u;\boldsymbol{\omega}),
\label{ell-3eq}\end{eqnarray} the latter of which can be rewritten
with the aid of modular transformation (\ref{mod}) as
\begin{equation}
G(u+\omega_3;\boldsymbol{\omega})=e^{-\pi iB_{2,2}(u;\mathbf{\omega})}
 G(u;\boldsymbol{\omega}),
\label{ell-3eq-r}\end{equation}
\end{subequations}
where
$$
B_{2,2}(u;\boldsymbol{\omega})=\frac{u^2}{\omega_1\omega_2}
-\frac{u}{\omega_1}-\frac{u}{\omega_2}+
\frac{\omega_1}{6\omega_2}+\frac{\omega_2}{6\omega_1}+\frac{1}{2}.
$$
Such a system
of three difference equations determines the meromorphic function
$G(u;\boldsymbol{\omega})$ up to a multiplicative constant (which
is the only meromorphic function with three pairwise
incommensurate periods $\omega_{1,2,3}$). Similar to the
$\Gamma(z;q,p)$ function, the $G(u;\boldsymbol{\omega})$ function
satisfies a simple reflection equation given by
\begin{equation}
G(u;\boldsymbol{\omega})G(-u;\boldsymbol{\omega})= \frac{e^{\pi
iB_{2,2}(u;\mathbf{\omega})}}
{\theta(e^{-2\pi i\frac{u}{\omega_1}};r)\theta(e^{-2\pi
i\frac{u}{\omega_2}};p)}.
\end{equation}

If we fix the quasiperiods  $\omega_1$, $\omega_2$ such that
$\text{Im}(\omega_{1}/\omega_2)>0$, and take $\omega_3$ to
infinity in such a way that $\text{Im}(\omega_{3}/\omega_1),
\text{Im}(\omega_{3}/\omega_2)\to +\infty$ (so $p,r\to 0$), then
we obtain
\begin{equation}
\lim_{p,r\to 0} \frac{1}{G(u;\boldsymbol{\omega})}
=S(u;\omega_1,\omega_2) = \frac{(e^{2\pi i u/\omega_2}; q)_\infty}
{(e^{2\pi iu/\omega_1}\tilde q; \tilde q)_\infty}.
\label{2d-sin}\end{equation}
In the modern time, the function $S(u;\boldsymbol{\omega})$ was
introduced by Shintani \cite{shi:kronecker}. It is related to the
Barnes double gamma function and is called the double sine function
\cite{kur:multiple}. Its various properties and applications are
discussed, e.g., in \cite{fad:modular,jim-miw:quantum,kls:unitary,
man:lectures,nis-uen:integral,rui:first,rui:generalized}. Some
$q$-beta integrals expressed in terms of
$S(u;\boldsymbol{\omega})$ were considered in
\cite{fkv:strongly,pon-tes:clebsch,sto:hyperbolic,tak:twisted}.

One of the main properties of the double sine function consists of
the fact that it can be extended continuously in the quasiperiods
$\omega_{1,2}$ from the upper half plane
$\text{Im}(\omega_{1}/\omega_2)>0$ (so $|q|<1$) to the positive
real axis $\omega_1/\omega_2>0$ (so $|q|=1$), the resulting
function still being meromorphic in $u$. A similar situation holds
for the function $G(u;\boldsymbol{\omega})$, as can be verified by
expressing it in terms of the Barnes' multiple gamma function of
the third order \cite{spi:integrals}. The following theorem
provides a convenient representation of $G(u;\boldsymbol{\omega})$
detailing explicitly its structure when $q$ lies on the unit
circle.

\begin{theorem}
Let $\text{Im}(\omega_3/\omega_1)$,
$\text{Im}(\omega_3/\omega_2)>0$ and
$\text{Im}(\omega_1/\omega_2)\geq 0$. The analytic continuation of
$G(u;\boldsymbol{\omega})$ \eqref{ell-3eq-r} from the domain
$\text{Im}(\omega_1/\omega_2)>0$ to the boundary
$\omega_1/\omega_2>0$ is given by the meromorphic function
\begin{subequations}
\begin{equation}
G(u;\boldsymbol{\omega})=e^{-\pi iP(u)}\Gamma(e^{-2\pi
i\frac{u}{\omega_3}}; \tilde r, \tilde p),
\label{gamma-tr}\end{equation} where $P(u)$ is the following
polynomial of the third degree
\begin{equation}
P(u)=\frac{1}{3\omega_1\omega_2\omega_3}
\left(u-\frac{1}{2}\sum_{n=1}^3\omega_n\right)
\left(u^2-u\sum_{n=1}^3\omega_n+\frac{\omega_1\omega_2\omega_3}{2}
\sum_{n=1}^3\frac{1}{\omega_n}\right). \label{p(u)}\end{equation}
\end{subequations}
\end{theorem}

\begin{proof}
Let us first assume temporarily that
$\text{Im}(\omega_1/\omega_2)>0$. We denote the right-hand side of
(\ref{gamma-tr}) by $f(u)$. It is easy to see that
\begin{equation}
\frac{f(u+\omega_1)}{f(u)}=e^{\pi i(P(u)-P(u+\omega_1))}
\theta(e^{-2\pi i\frac{u}{\omega_3}};\tilde p) =\theta(e^{2\pi
i\frac{u}{\omega_2}}; p), \label{eq-new}\end{equation}
as a consequence of the
modular transformation for theta functions in \eqref{mod}. The
function $f(u)$ thus satisfies equation (\ref{ell-1eq}). By
symmetry, it satisfies (\ref{ell-2eq}) as well. Analogously, we
have that $f(u+\omega_3)/f(u)=e^{\pi i(P(u)-P(u+\omega_3))}$,
which is seen to coincide with (\ref{ell-3eq-r}). The independence
of the quasiperiods $\omega_{1,2,3}$ over $\mathbb{Q}$ now implies
that $G(u;\boldsymbol{\omega})/f(u)$ must be constant. Its value
is equal to one, because for $u=(\omega_1+\omega_2+\omega_3)/2$ we
have that $G(u;\boldsymbol{\omega})=f(u)=1$. The theorem then
follows upon analytic continuation of the right-hand side of
\eqref{gamma-tr} to the positive real axis $\omega_1/\omega_2>0$.
\end{proof}

\begin{remark}
In \cite{fel-var:elliptic}, modular transformation properties of
the elliptic gamma function were investigated. For
$|q|,|p|,|r|<1$ the equality in (\ref{gamma-tr}) follows from one of these modular
transformations. Function $\Gamma(z;q,p)$ has a pointwise limit for some particular
values of $\omega_1/\omega_2\in X\subset\mathbb{R}_+$ \cite{fel-var:elliptic},
but it does not assume validity of (\ref{gamma-tr}) in this regime.
Our result consists thus of the observation that,
after an appropriate rewriting, this modular transformation
provides a representation for the elliptic gamma function that is
well defined for all $\omega_1/\omega_2>0$ \cite{spi:integrals}.
\end{remark}

Below we shall need functional equations satisfied by
$S(u;\boldsymbol{\omega})$
\begin{equation}
\frac{S(u+\omega_1;\boldsymbol{\omega})}{S(u;\boldsymbol{\omega})}=
\frac{1}{1-e^{2\pi i \frac{u}{\omega_2}}},\qquad
\frac{S(u+\omega_2;\boldsymbol{\omega})}{S(u;\boldsymbol{\omega})}=
\frac{1}{1-e^{2\pi i \frac{u}{\omega_1}}}
\label{func-eq}\end{equation}
and its asymptotics
for $u$ going to infinity. Let us fix the incommensurate
quasiperiods $\omega_1$ and $\omega_2$ such that
$\text{Im}(\omega_1/\omega_2)>0$ (so $|q|<1$). Then it follows
from the infinite product representation (\ref{2d-sin}) and its
modular inverse
\begin{equation}
S(u;\omega_1,\omega_2) = e^{-\pi iB_{2,2}(u;\boldsymbol{\omega})}
\frac{(e^{-2\pi iu/\omega_1}; \tilde q)_\infty} {(e^{-2\pi i
u/\omega_2}q; q)_\infty}
\end{equation}
that
\begin{subequations}
\begin{eqnarray}\label{asymp1}
&&
\lim_{\text{Im}(\frac{u}{\omega_1}),\text{Im}(\frac{u}{\omega_2})\to
+\infty}S(u;\boldsymbol{\omega})=1,
\\ &&
\lim_{\text{Im}(\frac{u}{\omega_1}),\text{Im}(\frac{u}{\omega_2})\to
-\infty} e^{\pi iB_{2,2}(u;\boldsymbol{\omega})}
S(u;\boldsymbol{\omega}) =1 . \label{asymp2}\end{eqnarray}
\end{subequations}
It turns out that the same asymptotics also holds for the boundary
domain $\omega_1/\omega_2>0$ (so $|q|=1$), as can be verified by
means of an integral representation for $S(u;\boldsymbol{\omega})$
\cite{kls:unitary, rui:first}.

\begin{remark}
If we take the limit $p, r\to 0$ in (\ref{gamma-tr}), then the
transition from $G(u;\boldsymbol{\omega})$ to the double sine
function simplifies to
\begin{eqnarray}\nonumber
&& \lim_{\text{Im}(\frac{\omega_3}{\omega_1}),
\text{Im}(\frac{\omega_3}{\omega_2})\to +\infty} \left( e^{-\pi
i\omega_3 \frac{2u-\omega_1-\omega_2}{12\omega_1\omega_2} }
\Gamma(e^{-2\pi i\frac{u}{\omega_3}};e^{-2\pi
i\frac{\omega_3}{\omega_1}}, e^{-2\pi
i\frac{\omega_3}{\omega_2}})\right)
\\ && \makebox[4em]{}
= e^{-\pi i\frac{3(2u-\omega_1-\omega_2)^2-\omega_1^2-\omega_2^2 }
{24\omega_1\omega_2} } S^{-1}(u;\omega_1,\omega_2).
\label{gamma-lim}\end{eqnarray} Such a limiting relation was first
derived in a different way and in a stronger sense by
Ruijsenaars \cite{rui:first}.
\end{remark}

\section{The elliptic beta integral on the unit circle}

We now turn to the elliptic beta integrals. Let $|q|,|p|<1$ and
let $t_n$, $n=0, \dots,4$, be five complex parameters satisfying
the inequalities $|t_n|<1$ and $ |pq|<|A|$, where
$A\equiv\prod_{n=0}^4t_n$. The elliptic beta integral of
\cite{spi:elliptic} states that
\begin{equation}\label{ell-int}
\frac{1}{2\pi i}\int_\mathbb{T}\frac{\prod_{n=0}^4\Gamma(t_n z^\pm
;q,p)} {\Gamma(z^{\pm 2},A z^\pm ;q,p)}\frac{dz}{z}=
\frac{2\prod_{0\leq n<m\leq 4} \Gamma(t_nt_m;q,p)}
{(q;q)_\infty(p;p)_\infty\prod_{n=0}^4\Gamma(At_n^{-1};q,p)},
\end{equation}
where $\mathbb{T}$ denotes the positively oriented unit circle.
Here we have employed the shorthand notations
\begin{eqnarray*}
&& \Gamma(z_1,\ldots,z_m;q,p)\equiv\prod_{l=1}^m
\Gamma(z_l;q,p),\\
&& \Gamma(tz^\pm;q,p)\equiv\Gamma(tz,tz^{-1};q,p),\quad \Gamma(z^{\pm
2};q,p)\equiv\Gamma(z^2,z^{-2};q,p).
\end{eqnarray*}
For
$p= 0$, the integration formula (\ref{ell-int}) amounts to an
integral explicitly constructed by Rahman in \cite{rah:integral}
through a specialization of the Nassrallah-Rahman integral
representation for a very-well-poised basic hypergeometric
$_8\varphi_7$ series.

The following theorem provides a modified elliptic beta integral
that---unlike \eqref{ell-int}---is well defined for $|q|=1$.

\begin{theorem}\label{circle-int:thm}
Let $\text{Im} (\omega_1/\omega_2)\geq 0$ and $\text{Im}
(\omega_3/\omega_1), \text{Im} (\omega_3/\omega_2)>0$, and
let $g_n$, $n=0,\ldots,4$, be five complex parameters subject to
the constraints
$$
\text{Im}(g_n/\omega_3)<0, \qquad
\text{Im}((\mathcal{A}-\omega_1-\omega_2)/\omega_3)>0,
$$
with $\mathcal{A}\equiv \sum_{n=0}^4g_n$. Then the following
integration formula holds
\begin{subequations}
\begin{equation}
\int_{-\omega_3/2}^{\omega_3/2} \frac{\prod_{n=0}^4 G(g_n\pm
u;\boldsymbol{\omega})} {G(\pm 2u, \mathcal{A}\pm
u;\boldsymbol{\omega})} \frac{du}{\omega_2} = \kappa\, \frac{
\prod_{0\leq n<m\leq 4}G(g_n+g_m;\boldsymbol{\omega})}
{\prod_{n=0}^4G(\mathcal{A}-g_n;\boldsymbol{\omega})},
\label{circle-int}\end{equation} where
\begin{equation}
\kappa= \frac{-2(\tilde q;\tilde q)_\infty}
{(q;q)_\infty(p;p)_\infty(r;r)_\infty}. \label{kappa}
\end{equation}
\end{subequations}
Here the integration is meant along the straight line segment
connecting $-\omega_3/2$ to $\omega_3/2$ and we employed the
shorthand notation $G(a\pm
b;\boldsymbol{\omega})\equiv G(a+b,a-b;\boldsymbol{\omega})$.
\end{theorem}

\begin{proof}
We start by substituting  relation (\ref{gamma-tr}) into the left-hand side
of (\ref{circle-int}). This yields
\begin{equation}
e^{\pi ia}\int_{-\omega_3/2}^{\omega_3/2}\frac{\prod_{n=0}^4
\Gamma(e^{-2\pi i\frac{g_n\pm u}{\omega_3}};\tilde r,\tilde p)}
{\Gamma(e^{\mp 4\pi i\frac{u}{\omega_3}}, e^{-2\pi
i\frac{\mathcal{A}\pm u}{\omega_3}};\tilde r,\tilde p)}
\frac{du}{\omega_2}, \label{cir-int}\end{equation} where
\begin{eqnarray}\nonumber
&&
a=\frac{2}{3\omega_1\omega_2\omega_3}\left(\mathcal{A}^3-\sum_{n=0}^4g_n^3
\right) -\frac{\sum_{m=1}^3\omega_m}{\omega_1\omega_2\omega_3}
\left(\mathcal{A}^2-\sum_{n=0}^4g_n^2 \right)
\\  && \makebox[4em]{}
+\frac{1}{2}\left(\sum_{m=1}^3\omega_m\right)
\left(\sum_{m=1}^3\frac{1}{\omega_m}\right).
\nonumber\end{eqnarray} The constraints on the parameters permit
to employ formula (\ref{ell-int}) with the substitutions
\begin{equation*}
z\to e^{\frac{2\pi i}{\omega_3}u}, \quad  t_n\to e^{-\frac{2\pi i
}{\omega_3}g_n},\quad p\to e^{-2\pi i\frac{\omega_1}{\omega_3}},\quad
q\to e^{-2\pi i\frac{\omega_2}{\omega_3}} ,
\end{equation*}
which yields for (\ref{cir-int})
\begin{eqnarray*}
 \frac{2\omega_3\omega_2^{-1} e^{\pi ia}}{(\tilde r; \tilde r)_\infty
(\tilde p;\tilde p)_\infty}\frac{\prod_{0\leq n<m\leq 4}
\Gamma(e^{-2\pi i\frac{g_n+g_m}{\omega_3}};\tilde r,\tilde p)}
{\prod_{n=0}^4\Gamma(e^{-2\pi
i\frac{\mathcal{A}-g_n}{\omega_3}};\tilde r,\tilde p)} && \\
= \kappa\, \frac{\prod_{0\leq n<m\leq 4}
G(g_n+g_m;\boldsymbol{\omega})} {\prod_{n=0}^4
G(\mathcal{A}-g_n;\boldsymbol{\omega})}, && \nonumber
\end{eqnarray*}
with
$$
\kappa= \frac{2\omega_3e^{\frac{\pi i}{12}(\sum_{m=1}^3\omega_m)
(\sum_{m=1}^3\omega_m^{-1}) }}
{\omega_2(\tilde r; \tilde r)_\infty(\tilde p;\tilde p)_\infty}.
$$
After applying modular transformation (\ref{ded}) to the
infinite products appearing in $\kappa$, we obtain
$$
\kappa= -2\sqrt{\frac{\omega_1}{i\omega_2}}
\frac{e^{\frac{\pi i}{12}(\frac{\omega_1}{\omega_2}
+\frac{\omega_2}{\omega_1}) }}{(r;r)_\infty (p;p)_\infty}.
$$
One more
application of (\ref{ded}) allows us to replace the exponential
function by a ratio of infinite products, which entails the
desired form of $\kappa$ in (\ref{kappa}).
\end{proof}

Let us now consider the formal limit $p,r\to 0$ of the integral
(\ref{circle-int}). To this end, we fix the quasiperiods
$\omega_{1,2}$ such that $\text{Im}(\omega_1/\omega_2)\geq 0$ and
$\text{Re}(\omega_1/\omega_2)> 0$, and we furthermore pick
$\omega_3=it\omega_2$ with $t>0$. For $t\to +\infty$ the integral
of Theorem \ref{circle-int:thm} then degenerates formally to
\begin{equation}
\int_{\mathbb{L}}\frac{S(\pm 2u, \mathcal{A}\pm
u;\boldsymbol{\omega})} {\prod_{n=0}^4 S(g_n\pm
u;\boldsymbol{\omega})}\frac{du}{\omega_2}= -2\frac{(\tilde
q;\tilde q)_\infty}{(q;q)_\infty}
\frac{\prod_{n=0}^4S(\mathcal{A}-g_n;\boldsymbol{\omega})} {
\prod_{0\leq n<m\leq 4}S(g_n+g_m;\boldsymbol{\omega})},
\label{red1}\end{equation}
where the integration is along the line
$\mathbb{L}\equiv i\omega_2\mathbb{R}$, and with parameters
subject to the constraints $\text{Re}(g_n/\omega_2)>0$ and
$\text{Re}((\mathcal{A}-\omega_1)/\omega_2)<1$.
This integral was deduced by a similar formal limit from
the standard elliptic beta integral \eqref{ell-int} and
rigorously proved by Stokman in \cite{sto:hyperbolic},
where it was referred to as the hyperbolic Nassrallah-Rahman
integral.

\begin{remark}
The generalized gamma function $[z;\tau]_\infty$ used in
\cite{sto:hyperbolic} coincides with the double sine function
(\ref{2d-sin}) upon setting $z=(\omega_2-u)/\omega_2$ and
$\tau=-\omega_2/\omega_1$.
\end{remark}

\begin{remark}
To further elucidate the intimate relation between the elliptic beta integrals
\eqref{ell-int} and  \eqref{circle-int}, let us recall that---according
to the general definition introduced in
\cite{spi:integrals}---a contour integral is called an {\em elliptic
hypergeometric integral} if, for some displacement $\omega_1$,
the ratio of its integrands
$\Delta(u+\omega_1)/\Delta(u)$ constitutes
an elliptic function of $u$ with periods $\omega_2, \omega_3$ (say).
Now, by the change of variables
$$
z=e^{\frac{2\pi i}{\omega_2}u}, \qquad t_n=e^{\frac{2\pi i}{\omega_2}g_n},
\qquad A=e^{\frac{2\pi i}{\omega_2}\mathcal{A}},
$$
the integral
(\ref{ell-int}) passes into the additive form
\begin{subequations}
\begin{equation}\label{ell-int-add}
\int_{-\omega_2/2}^{\omega_2/2}\Delta(u)\frac{du}{\omega_2}=
\frac{2\prod_{0\leq n<m\leq 4} \Gamma(e^{2\pi
i\frac{g_n+g_m}{\omega_2}};q,p)}
{(q;q)_\infty(p;p)_\infty\prod_{n=0}^4\Gamma(e^{2\pi
i\frac{\mathcal{A}-g_n}{\omega_2}};q,p)},
\end{equation}
with the integrand given by
\begin{equation}
\Delta(u)=\frac{\prod_{n=0}^4\Gamma(e^{2\pi i \frac{g_n\pm
u}{\omega_2}};q,p)} {\Gamma(e^{\pm 4\pi
i\frac{u}{\omega_2}},e^{2\pi i\frac{\mathcal{A}\pm u}{\omega_2}};q,p)},
\label{int}\end{equation}
\end{subequations}
where $\Gamma(e^{a\pm b};q,p)\equiv\Gamma(e^{a+b},e^{a-b};q,p)$.
We then have that
\begin{equation}
\frac{\Delta(u+\omega_1)}{\Delta(u)}= e^{2\pi
i\frac{\omega_1}{\omega_2}}\frac{\theta(e^{4\pi
i\frac{u+\omega_1}{\omega_2}};p)\theta(e^{2\pi
i\frac{u+\omega_1-\mathcal{A}}{\omega_2}};p)} {\theta(e^{4\pi
i\frac{u}{\omega_2}};p)\theta(e^{2\pi i\frac{u+\mathcal{A}}{\omega_2}};p)}
\prod_{n=0}^4\frac{\theta(e^{2\pi i\frac{u+g_n}{\omega_2}};p)}
{\theta(e^{2\pi i\frac{u+\omega_1-g_n}{\omega_2}};p)} , \label{int-eq}
\end{equation}
which is seen to be an elliptic function of $u$ with the periods
$\omega_2$ and $\omega_3$. With the aid of the difference equation
\eqref{ell-1eq} for the modified gamma function, it is not
difficult to verify that the integrand of modified elliptic
beta integral \eqref{circle-int} also provides a solution to
\eqref{int-eq}. Hence both integrals \eqref{ell-int} and
\eqref{circle-int} are elliptic hypergeometric integrals with
integrands satisfying \eqref{int-eq}. Whereas the solution of
\eqref{int-eq} originating from integral
\eqref{ell-int} is well-defined only for $|q|<1$ (or for $|q|>1$
upon performing the inversion $q\to q^{-1}$ in the elliptic gamma
function \cite{spi:integrals}), the modified integral
\eqref{circle-int} corresponds to a solution that extends
analytically from the regime $|q|<1$ to the unit circle $|q|=1$.
In this sense the modified integral  \eqref{circle-int} may be
seen as a $|q|=1$ analog of the original elliptic beta integral
\eqref{ell-int}.
\end{remark}

\section{Multiple integrals}

The following integral constitutes a multidimensional
generalization of the elliptic beta integral in \eqref{ell-int}
\begin{eqnarray}\nonumber
\lefteqn{\frac{1}{(2\pi i)^N}\int_{\mathbb{T}^N} \prod_{1\leq
j<k\leq N} \frac{\Gamma(tz_j^\pm z_k^\pm;q,p)} {\Gamma(z_j^\pm
z_k^\pm;q,p)}
\prod_{j=1}^N\frac{\prod_{n=0}^4\Gamma(t_nz_j^\pm;q,p)}
{\Gamma(z_j^{\pm 2}, B z_j^\pm;q,p)}\frac{dz_1}{z_1}
\cdots \frac{dz_N}{z_N} } && \nonumber \\
&& = \frac{2^N N!}{(p;p)_\infty^N(q;q)_\infty^N}
\prod_{j=1}^N\frac{\Gamma(t^j;q,p)}{\Gamma(t;q,p)}
\frac{\prod_{0\leq n<m\leq 4}\Gamma(t^{j-1}t_nt_m;q,p)}
{\prod_{n=0}^4\Gamma(t^{1-j}t_n^{-1}B;q,p)},
\label{ds}\end{eqnarray} where $B\equiv t^{2N-2}\prod_{n=0}^4t_n$,
$\Gamma(tz^\pm x^\pm;q,p)\equiv
\Gamma(tzx,tzx^{-1},tz^{-1}x,tz^{-1}x^{-1};q,p)$, and with
parameters subject to the constraints $|p|,|q|, |t|, |t_n|<1$ and
$|pq|<|B|$. This multiple elliptic beta integral was first
formulated as a conjecture in \cite{die-spi:elliptic}. Next, it was
shown in \cite{die-spi:selberg} that the conjecture in question
follows from a vanishing hypothesis for a related multiparameter
elliptic beta integral. Recently, a complete proof of the integral
in \eqref{ds} was found by Rains \cite{rai:privite}.

The following theorem provides the corresponding multidimensional
generalization of the modified elliptic beta integral in Theorem
\ref{circle-int:thm}.

\begin{theorem}\label{mcircle-int:thm}
Let $\text{Im} (\omega_1/\omega_2)\geq 0$ and $\text{Im}
(\omega_3/\omega_1), \text{Im} (\omega_3/\omega_2) >0$, and
let $g$, $g_n$, $n=0,\ldots,4$, be six complex parameters subject
to the constraints
$$
\text{Im}(g/\omega_3),\text{Im}(g_n/\omega_3)<0, \qquad
\text{Im}((\mathcal{B}-\omega_1-\omega_2)/\omega_3)>0,
$$
with $\mathcal{B}\equiv (2N-2)g+\sum_{n=0}^4g_n$. Then
\begin{eqnarray}\nonumber
\int_{-\frac{\omega_3}{2}}^{\frac{\omega_3}{2}}\cdots
\int_{-\frac{\omega_3}{2}}^{\frac{\omega_3}{2}} \prod_{1\leq
j<k\leq N} \frac{G(g\pm u_j \pm u_k;\boldsymbol{\omega})} {G(\pm
u_j \pm u_k;\boldsymbol{\omega})} \prod_{j=1}^N\frac{
\prod_{n=0}^4 G(g_n \pm u_j;\boldsymbol{\omega}) } { G(\pm 2u_j,
\mathcal{B}\pm u_j;\boldsymbol{\omega})}\frac{du_1}{\omega_2}
\cdots \frac{du_N}{\omega_2}  && \nonumber \\  = \kappa^N
N!\prod_{j=1}^N\frac{G(jg;\boldsymbol{\omega})}{G(g;\boldsymbol{\omega})}
\frac{\prod_{0\leq n<m\leq
4}G((j-1)g+g_n+g_m;\boldsymbol{\omega})}
{\prod_{n=0}^4G((1-j)g+\mathcal{B}-g_n;\boldsymbol{\omega})},
\makebox[2em]{} && \label{ds-circle}\end{eqnarray} with $\kappa$
given by (\ref{kappa}) and $G(c\pm a \pm b;\boldsymbol{\omega})\equiv
G(c+a+b,c+a-b,c-a+b,c-a-b;\boldsymbol{\omega})$.
\end{theorem}
\begin{proof}
The proof is analogous to that of Theorem \ref{circle-int:thm}.
Specifically, after substituting
(\ref{gamma-tr}) into the left-hand side of (\ref{ds-circle}) and application
of the multiple beta integral in
(\ref{ds}), we arrive at the integration formula stated in the theorem
upon expressing the resulting evaluation
in terms of the modified gamma function $G(u;\boldsymbol{\omega})$.
To infer the correctness of
the value of the proportionality constant $\kappa^N N!$,
one observes that the dependence on $g$ in the factors
originating from the exponential multipliers cancels out.
The expression for
the proportionality constant then follows from the fact that
for $g\to 0$ integral \eqref{ds-circle} reduces to the $N$-th power of
the elliptic beta integral \eqref{circle-int}.
\end{proof}

\begin{remark}
In \cite{die-spi:selberg,spi:integrals} various other types of
multiple elliptic beta integrals were formulated.
These can be extended to the unit circle $|q|=1$ in a similar fashion.
\end{remark}

For $p=0$, elliptic beta integral \eqref{ds} reduces to a
Gustafson's multiple integral corresponding to the Nassrallah-Rahman
type generalization of the Selberg integral \cite{gus:some}.
The following theorem provides a corresponding
multiple analog of the integral in \eqref{red1}.  The integration
formula in question can be obtained formally from the modified
elliptic beta integral \eqref{ds-circle} with
$\omega_1/\omega_2>0$ by taking the limit $p,r\to 0$ in the manner
explained below Theorem \ref{circle-int:thm}.

\begin{theorem}\label{nr:thm}
Let $\omega_1,\omega_2$ be quasiperiods such that
$\text{Im}(\omega_1/\omega_2)\geq 0$ and
$\text{Re}(\omega_1/\omega_2)>0$. Furthermore, let $g, g_n,$
$n=0,\ldots,4$, be parameters subject to the constraints
$\text{Re}(g/\omega_1),\text{Re}(g/\omega_2),
\text{Re}(g_n/\omega_2)>0$ and
$\text{Re}((\mathcal{B}-\omega_1)/\omega_2)<1$ (with $\mathcal{B}$
as in Theorem \ref{mcircle-int:thm}). Then
\begin{subequations}
\begin{equation}\label{q-ds-circle}
\int_{\mathbb{L}^N} \Delta(\mathbf{u};\mathbf{g})
\frac{du_1}{\omega_2} \cdots \frac{du_N}{\omega_2}
=\mathcal{N}(\mathbf{g}),
\end{equation}
where $\mathbb{L}=i\omega_2\mathbb{R}$,
\begin{equation}
\Delta(\mathbf{u};\mathbf{g})=\prod_{1\leq j<k\leq N} \frac{S(\pm
u_j \pm u_k;\boldsymbol{\omega})} {S(g\pm u_j \pm
u_k;\boldsymbol{\omega})} \prod_{j=1}^N\frac{S(\pm 2u_j,
\mathcal{B} \pm u_j;\boldsymbol{\omega})} {\prod_{n=0}^4 S(g_n \pm
u_j;\boldsymbol{\omega}) } \label{q-deg}\end{equation} and
\begin{equation}
\mathcal{N}(\mathbf{g}) = (-2)^N N!\frac{(\tilde q;\tilde
q)_\infty^N}
{(q;q)_\infty^N}\prod_{j=1}^N\frac{S(g;\boldsymbol{\omega})}
{S(jg;\boldsymbol{\omega})}
\frac{\prod_{n=0}^4S((1-j)g+\mathcal{B}-g_n;\boldsymbol{\omega})}
{\prod_{0\leq n<m\leq 4}S((j-1)g+g_n+g_m;\boldsymbol{\omega})}.
\label{N}\end{equation}
\end{subequations}
\end{theorem}

Through a parameter specialization, Gustafson's multiple integral
of the Nassral\-lah-Rahman type can be reduced to a multiple
Askey-Wilson integral first derived in \cite{gus:generalization}.
The corresponding degeneration of Theorem \ref{nr:thm} reads as
follows.

\begin{theorem}\label{aw:thm}
Let $\omega_1,\omega_2$ be quasiperiods such that
$\text{Im}(\omega_1/\omega_2)\geq 0$ and
$\text{Re}(\omega_1/\omega_2)> 0$, and let $g, g_n,$
$n=0,\ldots,3$, be parameters subject to the constraints
$\text{Re}(g/\omega_1),$ $\text{Re}(g/\omega_2),$
$\text{Re}(g_n/\omega_2)>0$ and
$\text{Re}((\mathcal{B}-\omega_2)/\omega_1)<1$ with
$\mathcal{B}\equiv (2N-2)g+\sum_{n=0}^3g_n$. Then
\begin{eqnarray}\nonumber
\lefteqn{ \int_{\mathbb{L}^N} \prod_{1\leq j<k\leq N} \frac{S(\pm
u_j \pm u_k;\boldsymbol{\omega})} {S(g\pm u_j \pm
u_k;\boldsymbol{\omega})} \prod_{j=1}^N\frac{S(\pm
2u_j;\boldsymbol{\omega})} {\prod_{n=0}^3 S(g_n \pm
u_j;\boldsymbol{\omega}) }\frac{du_1}{\omega_2} \cdots
\frac{du_N}{\omega_2} } && \nonumber \\ && = (-2)^N
N!\frac{(\tilde q;\tilde q)_\infty^N}
{(q;q)_\infty^N}\prod_{j=1}^N\frac{S(g;\boldsymbol{\omega})}
{S(jg;\boldsymbol{\omega})}
\frac{S((1-j)g+\mathcal{B};\boldsymbol{\omega})} {\prod_{0\leq
n<m\leq 3}S((j-1)g+g_n+g_m;\boldsymbol{\omega})}.
\label{q-ds-selberg}\end{eqnarray}
\end{theorem}

For $N=1$, the integral in Theorem \ref{aw:thm} reduces to
a single variable Askey-Wilson type integral
\begin{equation}
\int_{\mathbb{L}}\frac{S(\pm 2u;\boldsymbol{\omega})}
{\prod_{n=0}^3 S(g_n\pm
u;\boldsymbol{\omega})}\frac{du}{\omega_2}= -2\frac{(\tilde
q;\tilde q)_\infty}{(q;q)_\infty}
\frac{S(g_0+g_1+g_2+g_3;\boldsymbol{\omega})} { \prod_{0\leq
n<m\leq 3}S(g_n+g_m;\boldsymbol{\omega})}, \end{equation}
which was established by Ruijsenaars \cite{rui:int}
and Stokman \cite{sto:hyperbolic}.

Formally, the integral of Theorem \ref{aw:thm} follows from that
of Theorem \ref{nr:thm} with $\omega_1/\omega_2>0$ upon setting
$g_4\to g_4+i\omega_2t$ and performing the limit $t\to +\infty$.
However, it is not a simple matter to upgrade such formal limiting
relations between Theorem \ref{mcircle-int:thm} and Theorems
\ref{nr:thm}, \ref{aw:thm} to rigorous proofs of the latter
integration formulas. A direct proof of Theorems \ref{nr:thm} and
\ref{aw:thm}, modelled after Rains' proof of the multiple elliptic
beta integral \eqref{ds}, is given in Section \ref{sec6} below.
As it was communicated to us by Stokman after finishing this paper,
Theorem \ref{aw:thm} can be proved by a multivariable generalization of
the method of \cite{sto:hyperbolic} as well.

\section{Proof of Theorems \ref{nr:thm} and
\ref{aw:thm}}\label{sec6}

We first detail the proof of Theorem \ref{nr:thm} and then
indicate some modifications so as to incorporate Theorem
\ref{aw:thm}.

Let us for the moment assume that the quasiperiods
$\omega_1,\omega_2$ are incommensurate over $\mathbb{Q}$. The
double sine function $S(u;\boldsymbol{\omega})$ then has simple
zeros located at the points
$u=-\omega_1\mathbb{N}-\omega_2\mathbb{N}$ and simple poles at
$u=\omega_1(1+\mathbb{N})+\omega_2 (1+\mathbb{N})$.
Therefore, the integrand $\Delta(\mathbf{u};\mathbf{g})$ in
\eqref{q-ds-circle} has poles at the points
\begin{eqnarray}\nonumber
&& \pm u_j=-\mathcal{B}+\omega_1(1+\mathbb{N})+\omega_2(1+\mathbb{N}),
\; g_n+\omega_1\mathbb{N}+\omega_2\mathbb{N},\; n=0,\ldots,4,
\\ \nonumber && \makebox[2em]{}
g+ u_k+\omega_1\mathbb{N}+ \omega_2\mathbb{N},\;
g- u_k+\omega_1\mathbb{N}+ \omega_2\mathbb{N},\;
k=1,\ldots,N,\: k\neq j,
\label{weight-poles}\end{eqnarray}
where $j=1,\ldots, N$.

Combination with the asymptotics in Eqs. \eqref{asymp1}, \eqref{asymp2}
reveals that the quotients
$S(u;\boldsymbol{\omega})/S(g+u;\boldsymbol{\omega})$ and
$S(2u,\mathcal{B}+u;\boldsymbol{\omega})/\prod_{n=0}^4S(g_n+u)$
are smooth and bounded on the complex line
$\mathbb{L}=i\omega_2\mathbb{R}$. Indeed, for $u=ix\omega_2$,
$x\in\mathbb{R}$ we stay away from poles and we have that
\begin{equation*}
\frac{S(i\omega_2x;\boldsymbol{\omega})}{S(g+i\omega_2x;\boldsymbol{\omega})}=
\begin{cases}
O(1) &\text{for}\; x\to +\infty \\
O(e^{-2\pi xg/\omega_1})&\text{for}\;x\to -\infty
\end{cases}
\end{equation*}
and
\begin{equation*}
\frac{S(2i\omega_2x,\mathcal{B}+i\omega_2x;\boldsymbol{\omega})}
{\prod_{n=0}^4S(g_n+i\omega_2x;\boldsymbol{\omega})}=
\begin{cases}
O(1) &\text{for}\; x\to +\infty \\
O(e^{2\pi x
(2(N-1)g/\omega_1+1+\omega_2/\omega_1)})&\text{for}\;x\to -\infty
\end{cases}.
\end{equation*}
It thus follows that the integrand $\Delta(\mathbf{u};\mathbf{g})$
is smooth  and exponentially decaying at infinity on the
integration domain $\mathbb{L}^N$. Hence, the integral in Eq.
\eqref{q-ds-circle} converges.

To infer the validity of the integration formula we distinguish
three parameters $g_0,g_1,g_2$ and factorize the integrand as
$\Delta(\mathbf{u};\mathbf{g})=
\Delta_+(\mathbf{u})\Delta_-(\mathbf{u})$ with
\begin{equation}
\Delta_+(\mathbf{u})=\prod_{1\leq j<k\leq N} \frac{S(u_j \pm
u_k;\boldsymbol{\omega})} {S(g + u_j \pm u_k;\boldsymbol{\omega})}
\prod_{j=1}^N\frac{S(2u_j, \mathcal{B} - u_j,\omega_1+ \mathcal{C}
-u_j; \boldsymbol{\omega})}{ S(\omega_1+ \mathcal{C}
+u_j;\boldsymbol{\omega})\prod_{n=0}^4
S(g_n+u_j;\boldsymbol{\omega})}, \label{d-plus}\end{equation}
where $\mathcal{C}=(N-1)g+g_0+g_1+g_2$ and $\Delta_-(\mathbf{u})=
\Delta_+(-u_1,\ldots,-u_n)$. Similarly, we introduce the shifted
functions
\begin{eqnarray}\label{tilde-d}
\lefteqn{ \tilde \Delta_+(\mathbf{u})=\prod_{1\leq j<k\leq N}
\frac{S(u_j \pm u_k;\boldsymbol{\omega})} {S(g + u_j \pm
u_k;\boldsymbol{\omega})} }&&
\\ && \times
\prod_{j=1}^N\frac{S( 2u_j,
\mathcal{B}+\frac{\omega_1}{2}-u_j,\mathcal{C}+\frac{\omega_1}{2}-u_j;
\boldsymbol{\omega})}{S(\mathcal{C}+\frac{\omega_1}{2}+u_j,g_3-\frac{\omega_1}{2}+u_j,
g_4-\frac{\omega_1}{2}+u_j;\boldsymbol{\omega}) \prod_{n=0}^2
S(g_n+\frac{\omega_1}{2}+u_j; \boldsymbol{\omega})},
\nonumber\end{eqnarray} and
$\tilde\Delta_-(\mathbf{u})=\tilde\Delta_+(-u_1,\ldots,-u_n)$,
which provide a factorization of the integrand for a shifted set
of parameter values:
$$
\tilde \Delta_+(\mathbf{u})\tilde \Delta_-(\mathbf{u})
=\Delta(\mathbf{u};g,
g_0+\frac{\omega_1}{2},g_1+\frac{\omega_1}{2},
g_2+\frac{\omega_1}{2}, g_3-\frac{\omega_1}{2},
g_4-\frac{\omega_1}{2}).
$$
Let us for the moment assume that the parameters are such that the
shifted parameters (as well as the ones obtained after shifts
by $\pm\omega_2/2$) also belong to the domain indicated in the
theorem. We then have the following equality
\begin{eqnarray}\nonumber
&& \int_{\mathbb{L}^N} \tilde
\Delta_+(u_1+\frac{\omega_1}{2},\ldots,u_N+\frac{\omega_1}{2})
\Delta_-(\mathbf{u})\, du_1\cdots du_N
\\ && \makebox[4em]{}
=\int_{\mathbb{L}^N}\tilde \Delta_+(\mathbf{u})
\Delta_-(u_1-\frac{\omega_1}{2},\ldots, u_N-\frac{\omega_1}{2})\,
du_1\cdots du_N, \label{integral-shift}\end{eqnarray} which
follows by shifting the integration contours $\mathbb{L}$ on the
left-hand side by $-\omega_1/2$. Notice that such shifts are permitted
by the Cauchy theorem due to the absence of poles in the strip
between $\mathbb{L}$ and $\mathbb{L}-\omega_1/2$ combined with the
exponential decay at infinity. Indeed, the quotient
$S(u;\boldsymbol{\omega})/S(g+u;\boldsymbol{\omega})$ is
holomorphic on the strip $\{ u=s\omega_1+ix\omega_2\mid 0\leq
s\leq 1,\; -\infty <x<\infty\}$ and the quotient $S( 2u,
\mathcal{B}+u,\mathcal{C}+u; \boldsymbol{\omega})/\left(
S(\mathcal{C}+\frac{\omega_1}{2}+u,g_3-\frac{\omega_1}{2}+u,
g_4-\frac{\omega_1}{2}+u;\boldsymbol{\omega}) \prod_{n=0}^2
S(g_n+\frac{\omega_1}{2}+u; \boldsymbol{\omega})\right)$ is
holomorphic on the strip $\{ u=s\omega_1+ix\omega_2\mid 0\leq
s\leq 1/2,\; -\infty <x<\infty\}$. By performing sign flips of the
form $u_j\to -u_j$ in the integration variables and summing over
all  $2^N$ possible ways, we obtain from \eqref{integral-shift}
that
\begin{eqnarray}\nonumber
&& \int_{\mathbb{L}^N} \rho(\mathbf{u};\mathbf{g})
\Delta_+(\mathbf{u})\Delta_-(\mathbf{u})\, du_1\cdots du_N
\\ && \makebox[4em]{}
=\int_{\mathbb{L}^N} \tilde\rho(\mathbf{u};\mathbf{g}) \tilde
\Delta_+(\mathbf{u})\tilde\Delta_-(\mathbf{u})\, du_1\cdots du_N,
\label{int-sum}\end{eqnarray} with
\begin{eqnarray*}\label{rho}
&& \rho(\mathbf{u};\mathbf{g})=\sum_{\nu_\ell=\pm 1\atop
\ell=1,\ldots,N}
\frac{\tilde\Delta_+(\nu_1u_1+\frac{\omega_1}{2},\ldots,
\nu_nu_N+\frac{\omega_1}{2})}{\Delta_+(\nu_1 u_1,\ldots,\nu_N u_N)},  \\
&& \tilde\rho(\mathbf{u};\mathbf{g})=\sum_{\nu_\ell=\pm 1\atop
\ell=1,\ldots,N}
\frac{\Delta_-(\nu_1u_1-\frac{\omega_1}{2},\ldots, \nu_N
u_N-\frac{\omega_1}{2})} {\tilde\Delta_-(\nu_1 u_1,\ldots,\nu_Nu_N)}.
\label{t-rho}\end{eqnarray*}
Simplification of $ \rho(\mathbf{u};\mathbf{g})$ reveals that
\begin{eqnarray}
  \rho(\mathbf{u};\mathbf{g}) &=& \sum_{\nu_\ell=\pm 1\atop \ell=1,\ldots,N}
\prod_{1\leq j<k\leq N}\frac{1-tz_j^{\nu_j}z_k^{\nu_k}}
{1-z_j^{\nu_j}z_k^{\nu_k}} \prod_{j=1}^N
\frac{(1-t^{N-1}t_0t_1t_2z_j^{-\nu_j})\prod_{n=0}^2(1-t_nz_j^{\nu_j})}
{1-z_j^{2\nu_j}} \nonumber\\
&=& \prod_{j=1}^{N}(1-t^{j-1}t_0t_1)(1-t^{j-1}t_0t_2)
(1-t^{j-1}t_1t_2),
\label{rho-sum} \end{eqnarray}
with $t=e^{2\pi ig/\omega_2}$, $t_n=e^{2\pi ig_n/\omega_2}$,
$z_k=e^{2\pi iu_k/\omega_2}$. The summand in \eqref{rho-sum} is
invariant under permutations of $z_k$ and inversions $z_k\to z_k^{-1}$.
The product of $\rho(\mathbf{u};\mathbf{g})$ and the factor
$$
\prod_{1\leq j<k\leq N}\frac{(1-z_jz_k)(1-z_jz_k^{-1})}{z_j}
\prod_{j=1}^N\frac{1-z_j^2}{z_j}
$$
yields a Laurent polynomial in $z_j,j=1,\ldots,N,$ which is
antisymmetric with respect to both transformations (separate permutations
of $z_j$ and inversions $z_j\to z_j^{-1}$). Any such polynomial is
proportional to the multiplicative factor given above. The constant
of proportionality is found after setting $z_j=t_0t^{N-j}$, which
leaves only one term in the sum (with all $\nu_j=1$) giving the
right-hand side expression in \eqref{rho-sum}.

In the same way, one obtains for
$\tilde\rho(\mathbf{u};\mathbf{g})$ that replacing in relation
\eqref{rho} $t_{0,1,2}$ by $t_{3,4}q^{-1/2}, t^{N-1}t_0t_1t_2q^{1/2}$,
\begin{eqnarray*}
\tilde\rho(\mathbf{u};\mathbf{g})&=&  \sum_{\nu_\ell=\pm 1\atop
\ell=1,\ldots,N} \prod_{1\leq j<k\leq
N}\frac{1-tz_j^{\nu_j}z_k^{\nu_k}} {1-z_j^{\nu_j}z_k^{\nu_k}}
\prod_{j=1}^N \Biggl(
\frac{(1-t_3q^{-1/2}z_j^{\nu_j})(1-t_4q^{-1/2}z_j^{\nu_j})}
{1-z_j^{2\nu_j}} \nonumber \\ && \makebox[6em]{}\times
(1-t^{N-1}t_0t_1t_2q^{1/2}z_j^{\nu_j})(1-Bq^{-1/2}z_j^{-\nu_j})\Biggr)
 \nonumber \\
&=& \prod_{j=1}^{N}(1-t^{j-1}t_3t_4/q)
(1-t^{1-j}B/t_3)(1-t^{1-j}B/t_4),
\end{eqnarray*}
where $B=e^{2\pi i\mathcal{B}/\omega_2}$.
The derived expressions demonstrate that the functions in question
are constant in the integration variables $\mathbf{u}$, that is
$\rho(\mathbf{u};\mathbf{g})=\rho(\mathbf{g})$ and
$\tilde{\rho}(\mathbf{u};\mathbf{g})=\tilde{\rho}(\mathbf{g})$.
Hence, we can pull the corresponding factors out of the integrals
and rewrite \eqref{int-sum} as
\begin{equation*}
\int_{\mathbb{L}^N} \Delta_+(\mathbf{u})\Delta_-(\mathbf{u})\,
du_1\cdots du_N =\frac{\tilde\rho(\mathbf{g})}{\rho(\mathbf{g})}
\int_{\mathbb{L}^N} \tilde
\Delta_+(\mathbf{u})\tilde\Delta_-(\mathbf{u})\, du_1\cdots du_N,
\end{equation*} whence
\begin{equation}\label{int-trans}
\frac{\mathcal{N}(\mathbf{g})}
{\mathcal{N}(g,g_0+\frac{\omega_1}{2},g_1+\frac{\omega_1}{2},
g_2+\frac{\omega_1}{2},g_3-\frac{\omega_1}{2},g_4-\frac{\omega_1}{2})}=
\frac{\tilde\rho(\mathbf{g})}{\rho(\mathbf{g})}.
\end{equation}

As a result, we deduce that the ratio of the left- and right-hand sides
of \eqref{q-ds-circle} is invariant with respect to the shifts
$g_{0,1,2}\to g_{0,1,2}+\omega_1/2, g_{3,4}\to
g_{3,4}-\omega_1/2$, and, by symmetry, any permutation of indices
of the parameters. The double sine function $S(u;\boldsymbol{\omega})$
and, so, the integrand in \eqref{q-ds-circle}
and the integral's value $\mathcal{N}(\mathbf{g})$ are symmetric
with respect to the permutation of $\omega_1$ and $\omega_2$
\cite{kls:unitary}. The contour of integration $\mathbb{L}$
breaks the symmetry between $\omega_{1,2}$, but the transformations
used in \eqref{d-plus}-\eqref{int-trans} are purely algebraic
and do not depend on the contour of integration. Therefore, the ratio
of interest is invariant with respect to the shifts
$g_{0,1,2}\to g_{0,1,2}+\omega_2/2,$ $g_{3,4}\to g_{3,4}-\omega_2/2$
as well (with all permutations of indices of parameters).

By analyticity, without changing the integral's value
we can replace the contour of integration $\mathbb{L}$ by any other
contour embracing the same set of poles. For an appropriately
deformed contour, we can establish invariance of
the ratio of interest under the shifts $g_n\to g_n+ k\omega_1/2+
m\omega_2/2$, for arbitrary $k,m\in\mathbb{Z}$.
Taking $\omega_1,\omega_2 > 0$, we can choose a subset of these points
with a limiting point in the parameter space for which we can choose
$\mathbb{L}$ as the integration contour. Moreover, making
sequential $\omega_{1,2}/2$ shifts in different directions we can escape
large intermediate deformations of the integration contour
(such an argument is similar to the one given in \cite{spi:elliptic}).
Therefore, the ratio of the left- and
right-hand sides of \eqref{q-ds-circle} is equal to a
function of $\omega_{1,2}$ and $g$, which we denote as
$f(\omega_1,\omega_2,g)$.

In order to see that $f(\omega_1,\omega_2,g)$ actually equals to
one, it is necessary to use an analog of the residue formula
derived in \cite{die-spi:elliptic}. Namely, we take one of the
parameters, say, $g_0$ such that one pole from the half plane
$\text{Re}(u/\omega_2)>0$ crosses the contour of integration
$\mathbb{L}$. A similar move takes place for the pole at $u=-g_0$
from the $\text{Re}(u/\omega_2)<0$ half plane
(due to the reflection invariance). We keep intact all other poles in
the half planes to the left or right of $\mathbb{L}$. This
is possible to do by an appropriate choice of $g, g_1,\ldots,
g_4$, and $\omega_{1,2}$. It is not difficult to see that the
residues of these crossing poles taken, say, over the variable
$u_N$ have poles at the points $u_k=\pm(g_0+g), k=1,\ldots,N-1,$
(instead of $u_k=\pm g_0$) and they are still located to the left
or right of $\mathbb{L}$ for $\text{Re}((g_0+g)/\omega_2)<0$. Similar
shifts of the poles occur each time we calculate the residues.
Therefore we denote $\rho_k=g_0+(k-1)g$ and impose the
restrictions
$$
\text{Re}\left(\frac{\rho_k}{\omega_2}\right)<0,
\quad k=1, \ldots, N,\qquad
\text{Re}\left(\frac{g_0+\omega_1}{\omega_2}\right),
\text{Re}\left(\frac{g_0+\omega_2}{\omega_2}\right)>0.
$$

Because of the taken constraints upon $g, \omega_{1,2}$, we get a
simpler residue formula than the one derived in
\cite{die-spi:elliptic}:
\begin{eqnarray}\nonumber
\lefteqn{
\int_{\mathbb{L}^N_d } \Delta(\mathbf{u};\mathbf{g})
\frac{du_1}{\omega_2}\cdots\frac{du_N}{\omega_2} } && \\ &&
=\sum_{m=0}^N 2^mm! \binom{N}{m} \int_{\mathbb{L}^{N-m}}
\mu_m(\mathbf{u}) \frac{du_1}{\omega_2}\cdots\frac{du_{N-m}}{\omega_2},
\label{res}\end{eqnarray}
where $\mathbb{L}_d$ is a deformation of the
contour $\mathbb{L}$ such that it separates the same sets of poles
as $\mathbb{L}$ did before we started to change $g_0$. The factor
$2^m$ emerges because the residues appear in pairs and their values
coincide (due to the $u_k\to -u_k$ reflection invariance of the integrand
and different orientation of the contours encircling poles to the
left and right of $\mathbb{L}$). The factors $m!$ and $\binom{N}{m}$
count the number of orderings of $m$ cycles and the number of ways
to pick up these cycles out of $N$ possibilities.

The residue functions have the form $\mu_0(\mathbf{u})=\Delta(\mathbf{u};\mathbf{g})$
and for $m>0$
\begin{equation}
\mu_m(\mathbf{u}) = \kappa_m
\delta_{m,N-m}(\mathbf{u})\Delta_{N-m}(\mathbf{u};\mathbf{g}),
\label{mu}\end{equation} where
$\Delta_{N-m}(\mathbf{u};\mathbf{g})$ is obtained from integrand
\eqref{q-deg} if we replace in it $N$ by $N-m$ but keep
$\mathcal{B}=(2N-2)g+\sum_{n=0}^4 g_n$ unchanged. Other
coefficients are
\begin{equation*}
\kappa_m =(-1)^m\frac{(\tilde q;\tilde
q)_\infty^m}{(q;q)_\infty^m} \prod_{1\leq j<k\leq m}
\frac{S(\pm\rho_k - \rho_j;\boldsymbol{\omega})}
{S(g\pm\rho_k-\rho_j;\boldsymbol{\omega})} \prod_{l=1}^m
\frac{S(-2\rho_l,\mathcal{B}\pm \rho_l;\boldsymbol{\omega})}
       {\prod_{n=1}^4S(g_n\pm\rho_l;\boldsymbol{\omega})},
\end{equation*}
and
\begin{equation*}
\delta_{m,N-m}(\mathbf{u}) = \prod_{\begin{subarray}{c}1\leq j\leq
m\\ 1\leq k\leq N-m\end{subarray}} \frac{S(\pm \rho_j\pm
u_k;\boldsymbol{\omega})} {S(g\pm\rho_j\pm
u_k;\boldsymbol{\omega})}.
\end{equation*}
The expressions for $\mu_m(\mathbf{u})$ are derived by induction.
Indeed, the form of $\mu_1(\mathbf{u})$ is easily established
after taking into account the relation
$$
\lim_{u\to \pm g_0}\frac{u\mp g_0}{S(g_0\mp
u;\boldsymbol{\omega})} =\pm \frac{\omega_2}{2\pi i}\frac{(\tilde
q;\tilde q)_\infty} {(q;q)_\infty}
$$
and the fact that the contours encircling the corresponding poles
are oriented clockwise for the upper signs and anticlockwise for
the lower signs (this gives the total minus sign in $\kappa_1$).

Suppose that $\mu_m$ is given by \eqref{mu} for some $m>1$. In
order to find $\mu_{m+1}$ it is necessary to compute the residues
for poles located at $u_{N-m}=\pm \rho_{m+1}$. A simple
computation shows that, indeed,
$$
\int_{c_m}
\mu_m(\mathbf{u})\frac{du_{N-m}}{\omega_2}=\mu_{m+1}(\mathbf{u}),
$$
where $c_m$ is a small size clockwise orientated contour
encircling the pole at $u_{N-m}=\rho_{m+1}$.

By analyticity, our deformations of the parameters and of the
contour of integration do not change the integral value and,
therefore, the right-hand side sum in \eqref{res} equals to
$f(\omega_1,\omega_2,g)\mathcal{N}(\mathbf{g})$. We now divide
both sides of this equality by $\mathcal{N}(\mathbf{g})$ and take
the limit $g_4\to -g_0-(N-1)g$. For $m<N$, the coefficients
$\kappa_m(\mathbf{g})$, which can be represented in the form
\begin{eqnarray*}
\lefteqn{ \kappa_m =(-1)^m\frac{(\tilde q;\tilde q)_\infty^m}
{(q;q)_\infty^m}\prod_{l=1}^m \Biggl(
\frac{S(g,(2-m-l)g-2g_0;\boldsymbol{\omega})}{S(lg;\boldsymbol{\omega})}
}
\\ && \times
\frac{S(2g_0+\sum_{r=1}^4g_r+(2n+l-3)g,
\sum_{r=1}^4g_r+(2n-l-1)g;\boldsymbol{\omega})}
{\prod_{r=1}^4S(g_r+g_0+(l-1)g,g_r-g_0-(l-1)g;
\boldsymbol{\omega})} \Biggr),
\end{eqnarray*}
do not contain diverging factors in this limit and the integrals,
which they are multiplied by, remain bounded. Therefore, only the
term with $m=N$ survives and, by simple computation, we obtain
$$
\lim_{g_4\to
-g_0-(N-1)g}\frac{2^NN!\kappa_N}{\mathcal{N}(\mathbf{g})}=1,
$$
which means that $f(\omega_1,\omega_2,g)=1$. After proving
equality \eqref{q-ds-circle} in the taken restricted region of
parameters (where the parameters shifted by $\pm\omega_{1,2}/2$
satisfy the needed
constraints and $\omega_{1,2}>0$), we can analytically extend
it to the values of $\omega_{1,2}$ and parameters $g, g_n$
in the domain indicated in the formulation of the theorem.
Theorem \ref{nr:thm} is thus proved.

We now turn to the Askey-Wilson type integral \eqref{q-ds-circle}.
Its convergence conditions essentially differ from the previous
case. Indeed, we have
\begin{equation*}
\frac{S(2i\omega_2x;\boldsymbol{\omega})}
{\prod_{n=0}^3S(g_n+i\omega_2x;\boldsymbol{\omega})}=
\begin{cases}
O(1) &\text{for}\; x\to +\infty \\
O(e^{2\pi x(1+\omega_2/\omega_1-
\sum_{n=0}^3 g_n/\omega_1 )}) &\text{for}\;x\to -\infty
\end{cases}.
\end{equation*}
Combining together these limiting relations with the asymptotics for
the ratio
$S(i\omega_2x;\boldsymbol{\omega})/S(g+i\omega_2x;\boldsymbol{\omega})$,
we see that the integrand remains bounded in the integration domain
$\mathbb{L}^N$ and decays exponentially fast on its infinities if we
take $\text{Re}((\mathcal{B}-\omega_2)/\omega_1)<1$, where
$\mathcal{B}=(2N-2)g+\sum_{n=0}^3g_n$.

Invariance of the ratio of left- and right-hand sides of equality
(\ref{q-ds-circle}) under the specified parameter shifts relied only on
algebraic manipulations with the integrand. Therefore we can repeat them
for the limiting expression of the integrand appearing after
taking the limit $\text{Im}(g_4/\omega_2),$ $\text{Im}(g_4/\omega_1)\to
+\infty$ (or $t_4\to 0$). This simplifies the integrand
for (\ref{q-ds-circle}) to the one for \eqref{q-ds-selberg}.
Therefore, limiting analogs of equalities
\eqref{d-plus}--\eqref{int-trans} show that the ratio of the left-
and right-hand sides of \eqref{q-ds-selberg} do not depend on the
shifts in the parameter space $g_{0,1,2}\to g_{0,1,2}+\omega_{1,2}/2,$
$g_3\to g_3-\omega_{1,2}/2$ and the ones obtained by permutation of
indices. Using, again, an analytical continuation and the appropriately
simplified version of the residue calculus, we see that the ratio of
interest is actually equal to one. As a result, we establish validity
of integral (\ref{q-ds-selberg}) as well.

\medskip
\centerline{Acknowledgments}
\smallskip

The second author thanks the Instituto de Matem\'atica y
F\'{\i}sica of the Universidad de Talca for the hospitality during
the visit in May 2003, at which time the main results of this paper
were obtained. This work is supported in part by the Fondo
Nacional de Desarrollo Cient\'{\i}fico y Tecnol\'ogico (FONDECYT)
Grants No. \# 1010217 and No. \# 7010217, by the Programa Formas
Cuadr\'aticas of the Universidad de Talca, and by the Russian
Foundation for Basic Research (RFBR) Grant No. 03-01-00781.

\bibliographystyle{amsplain}

\end{document}